\documentclass[12pt]{amsart}   

\usepackage{amsmath, amssymb}   
\usepackage{hyperref}
\usepackage{url}
\newtheorem{thm}{Theorem}[subsection]
\newtheorem{prop}[thm]{Proposition}
\newtheorem{lem}[thm]{Lemma}
\newtheorem{conj}[thm]{Conjecture}
\newtheorem{cor}[thm]{Corollary}

\def\E{\mathbb{E}}
\def\P{\mathbb{P}}
\def\V{\mathbb{V}}
\def\Po{\mathcal{P}}

\usepackage{fullpage}

\usepackage{xypic}

\begin{document}

\title[Compositions of random involutions]{The cycle structure of compositions of random involutions}

\date{November 18, 2009}

\author{Michael Lugo}

\address{Department of Mathematics, University of Pennsylvania, 209 South 33rd Street, Philadelphia, PA 19104}

\email{mlugo@math.upenn.edu}



\maketitle

\begin{abstract}
In this article we consider the cycle structure of compositions of pairs of involutions in the symmetric group $S_n$ chosen uniformly at random.  These can be modeled as modified $2$-regular graphs, giving rise to exponential generating functions.  A composition of two random involutions in $S_n$ typically has about $\sqrt{n}$ cycles, and the cycles are characteristically of length $\sqrt{n}$.  Compositions of two random fixed-point-free involutions, on the other hand, typically have about $\log n$ cycles and are closely related to permutations with all cycle lengths even.  The number of factorizations of a random permutation into two involutions appears to be asymptotically lognormally distributed, which we prove for a closely related probabilistic model.  This study is motivated by the observation that the number of involutions in $[n]$ is $\sqrt{n!}$ times a subexponential factor; more generally the number of permutations with all cycle lengths in a finite set $S$ is $n!^{1-1/m}$ times a subexponential factor, and the typical number of $k$-cycles is nearly $n^{k/m}/k$.  Connections to pattern avoidance in involutions are also considered.
\end{abstract}

  \subsection{Introduction}

The purpose of this paper is to study the cycle structure of compositions of involutions.  Recall that an {\it involution} is a permutation with all cycles having length $1$ or $2$.  Let $a_n$ be the number of involutions in the symmetric group $S_n$.  Then as $n \to \infty$, 
\begin{equation}\label{eq:number-of-involutions} a_n \sim \sqrt{n!} e^{\sqrt{n}} (8 \pi e n)^{-1/4}. \end{equation}
This form involving $\sqrt{n!}$ is due to \cite[p. 583]{flajolet-sedgewick}; see \cite{moser-wyman} for the result in another form, and \cite[Example 3.2]{pemantle-notes} for details of the asymptotic analysis by the saddle-point method.  The factor $\sqrt{n!}$ is much faster-growing than $e^{\sqrt{n}} (8\pi en)^{-1/4}$.  So in a logarithmic sense the number of involutions of $[n]$ is approximately the square root of the number of permutations of $[n]$.  Thus the number of pairs of involutions of $[n]$ is logarithmically near $n!$.  This suggests identifying permutations with pairs of involutions.  A natural way to combine two involutions to form a permutation is composition, so we study compositions.  We also seek other ways in which involutions and permutations are related by forming ordered pairs.  In particular, Stanley-Wilf limits for various classes of pattern-avoiding permutations are known; in those cases where a Stanley-Wilf limit for the corresponding pattern-avoiding involutions exists, the former is the square of the latter.

We begin by asymptotically enumerating permutations with all cycle lengths in a finite set $S$; involutions are the case $S = \{1, 2 \}$.  Call a permutation with all cycle lengths in $S$ an $S$-permutation.  Let $p_n^{(S)}$ be the probability that a permutation of $[n]$ selected uniformly at random is an $S$-permutation.  Then $\lim_{n \to \infty} {\log p_n^{(S)} \over \log n!} = -1/(\max S)$; a refinement of this is Theorem \ref{thm:leading-term} below.  In particular the number of $k$-cycles of a typical $S$-permutation is near ${1 \over k} n^{k/(\max S)}$.  A typical involution of $[n]$ has $\sqrt{n}$ fixed points.  We then proceed to represent involutions graph-theoretically as partial matchings; thus compositions of two involutions can be identified with graphs having $2$-colored edges, where each vertex has at most one incident edge of each color.  The components of such graphs are paths and cycles, so we easily find generating functions involving them. This is our principal tool for extracting information on the cycle structure of these graphs and the corresponding permutations.  In particular, if $\sigma$ and $\tau$ are random involutions of $[n]$, then as $n \to \infty$:

\begin{itemize}
 \item The distribution of the number of $k$-cycles of $\tau \circ \sigma$ converges in distribution to $\Po(1) + 2\Po(1/2k)$ (Theorem \ref{thm:k-cycle-dist});
\item The mean number of cycles of $\tau \circ \sigma$ (of all lengths) is $\sqrt{n} + {1 \over 2} \log n + O(1)$ (Theorem \ref{thm:mean-number-permutation-cycles});
\item If $\sigma$ and $\tau$ are constrained to be fixed-point-free, then the distribution of the number of cycles of $\tau \circ \sigma$ is asymptotically normal with mean $\log n$ and variance $2 \log n$.
\end{itemize}

Finally, we consider the number of factorizations of a permutation into involutions.  The mean number of factorizations of $\pi$ into two involutions -- that is, solutions to $\pi = \tau \circ \sigma$, with $\tau$ and $\sigma$ involutions of $[n]$ -- is $e^{2\sqrt{n}}/\sqrt{8\pi e n} (1 + o(1))$.  We derive a formula (Theorem \ref{thm:inv-factorizations}) for the number of factorizations of $\pi \in S_n$ into two involutions, in terms of the cycle type of $\pi$.  This is a product over cycle lengths.  In a model of random permutations in which there are $\Po(1/k)$ cycles of length $k$ for $k = 1, \ldots, n$, the number of factorizations of a random permutation $\pi$ is lognormally distributed (Theorem \ref{thm:lognormal}).  If $\P_n^*$ denotes this probability measure, $F(\pi)$ the number of factorizations of $\pi$, and $\Phi$ the standard normal cdf, then
\[ \lim_{n \to \infty} \P_n^* \left( {\log F(\pi) - {1 \over 2} (\log n)^2 \over {1 \over 3} (\log n)^3} \le x \right) \to \Phi(x). \]
In particular the median number of factorizations of $\pi$ is near $\exp((\log n)^2/2)$, much smaller than the mean.  This is one of many indications that the measure on $S_n$ coming from compositions of involutions chosen uniformly at random is much different from the uniform measure on $S_n$.

\subsection{Pattern avoidance.} We recall the Stanley-Wilf conjecture (now the Marcus-Tardos theorem \cite{marcus-tardos}) on pattern avoidance.  Let $S_n(\pi)$ denote the set of $\pi$-avoiding permutations of $[n]$; then $S_n(\pi)$ is bounded above by $C^n$, for some constant $C$ depending on $\pi$.  We call the smallest such $C$ the {\it growth rate} of the pattern $\pi$ and denote it by $L(\pi)$.  Arratia \cite{arratia} has shown that the Marcus-Tardos theorem is equivalent to the existence of the limit $\lim_{n \to \infty} |S_n(\pi)|^{1/n}$, which equals $L(\pi)$, for all patterns $\pi$.  
  Now let $I_n(\pi)$ denote the set of $\pi$-avoiding {\it involutions} of $[n]$.  Then we can define the {\it involutory growth rate} of a pattern, $L_i(\pi) = \lim_{n \to \infty} |I_n(\pi)|^{1/n}$.  This limit may not exist in general, but it does in some special cases, leading to the following conjecture. \begin{conj}\label{conj:involution-growth} Let $\pi$ be a permutation pattern.  Then $L_i(\pi)$ exists and $L_i(\pi)^2 = L(\pi)$.   \end{conj}
\begin{table}
 \begin{tabular}{p{2cm}|p{4cm}|p{4cm}|p{3cm}}
$\pi$ & $I_n(\pi)$ & $S_n(\pi)$ & $I_n(\pi)^2/S_n(\pi)$ \\
\hline
$12 \ldots k$ & $\sim a_k (k-1)^n \newline (1/n)^{(k-1)(k-2)/4}$ 
 & $\sim b_k (k-1)^{2n} \newline (1/n)^{k^2/2-k}$ 
 & $\sim c_k n^{-1+k/2}$ \\
 & \cite[4.5 Case 1]{regev}, \cite{gessel} & \cite[4.5 Case 2]{regev} & \\
\hline
1234, 2143, 3412, 4321, 1243 & $M_n \sim \sqrt{27 \over 4\pi} 3^n n^{-3/2}$ & $\sim {81\sqrt{3} \over 16\pi} 9^n n^{-4}$ & ${4 \sqrt{3} \over 9} n$ \\
& \cite{egge-mansour} & \cite[Cor. 3.1.7]{west} shows patterns are Wilf-equivalent & \\
\hline
123, 132, 213, 321 & ${n \choose \lfloor n/2 \rfloor} \sim 2^n/\sqrt{\pi n}$ & $C_n \sim 4^n/\sqrt{\pi n^3}$ & $\sqrt{n/\pi}$ \\
 & \cite{simion-schmidt} & & \\
\hline
231, 312 & $2^{n-1}$ \cite{simion-schmidt} & $C_n \sim 4^n/\sqrt{\pi n^3}$ & $\sqrt{\pi \over 16} n^{3/2}$ \\
\hline
54321 & $C_{\lceil n/2 \rceil} C_{1 + \lfloor n/2 \rfloor} \sim {32 \over \pi} {4^n \over n^3}$ & $2^{25/2} 3 \pi^{-3/2} 16^n n^{-15/2}$ & ${1 \over 24} \sqrt{2 \over \pi} n^{3/2}$ \\
 & \cite{bousquet-melou} & \cite{regev} and symmetry & \\
\hline
\end{tabular}
\caption{Table of patterns for which the ordinary and involutory growth rates are both known. $C_n$ and $M_n$ are the Catalan and Motzkin numbers, respectively; $a_k, b_k, c_k$ are constants depending on $k$.}\label{table:involutions}
\end{table}
Table \ref{table:involutions} shows $I_n(\pi)$, $S_n(\pi)$, and the ratio of their squares in cases when both are known.   We note in particular that the conjecture is true for all patterns of length at most 3.  We also note that Wilf-equivalence of two patterns is not the same as ``involutory Wilf-equivalence''.  In particular $|S_n(\pi^r)| = |S_n(\pi)|$, where $\pi^r$ is the reversal of $\pi$, but it is not necessarily true that $|I_n(\pi^r)| = |I_n(\pi)|$.  Counterexamples include $\pi = 132$ and $\pi = 12345$; we have $|I_n(12345)| \sim (\pi^3/8) 4^n n^{-3}$ \cite{regev} but $I_n(54321) \sim {32 \over \pi} 4^n n^{-3}$ \cite{bousquet-melou}.  

The pattern $1342$ has growth rate 8 and the pattern $12453$ has growth rate $(1+\sqrt{8})^2$ \cite{bona}; the latter is the first known example of a pattern with non-integer growth rate.  Bona has shown \cite[Lemma 5.4]{bona} that given a pattern $\pi$ of growth rate $L(\pi) = g^2$, the pattern $\pi^\prime$ obtained by adding 1 to each element of $\pi$ and prepending 1 to the result has growth rate $L(\pi^\prime) = (g+1)^2$.   In other words, this operation raises the square root of the growth rate by 1; thus there is some precedent for studying $\sqrt{L(\pi)}$.  Perhaps in general $L_i(\pi^\prime) = L_i(\pi) + 1$.

Conjecture \ref{conj:involution-growth} can be restated probabilistically.  The probability that a random permutation of $[n]$ is $\pi$-avoiding seems to be the square of the probability that a random involution of $[n]$ is $\pi$-avoiding, multiplied by some asymptotically subexponential factor.  (In the few known cases this factor is $Cn^{-k}$ for some real constant $C$ and nonnegative rational number $k$.)  Thus involutions are, in general, more likely to avoid patterns than ordinary permutations.  This is because an involution is, in a sense, half a permutation.  The RSK algorithm \cite{stanley-99} takes a permutation $\pi$ to a pair of Young tableaux $(P, Q)$; if $\pi$ is an involution then $P = Q$, so involutions can be identified with individual Young tableaux.  The ``graph'' of a permutation $\pi$ is the set of points $\{ (i, \pi(i)) : 1 \le i \le n \}$ and an involution can be specified by fixing only the points on or below the diagonal, identifying involutions with half-graphs.  

Finally, Egge has studied permutations with graphs which are symmetric under other reflections or rotations \cite{egge}.  One might hope these lead to further generalizations of Conjecture \ref{conj:involution-growth}.  To give an example, involutions invariant under the reverse complement are determined by one-fourth of their graph, and the number of such permutations which are also $132$-avoiding grows like $2^{n/2}$.  Up to polynomial factors this is the fourth root of the Catalan number $C_n$, which is the number of $132$-avoiding permutations.  
Wulcan \cite{wulcan} has enumerated involutions avoiding generalized patterns, including all the generalized patterns of length 3; at this point no systematic review of the growth rates of the corresponding patterns in permutations has been undertaken.

\subsection{The number of permutations with all cycle lengths in some finite set}

The fact that the number of involutions of $[n]$ is approximately $\sqrt{n!}$ can be generalized to permutations with cycle lengths lying in any finite set.    We call a permutation with all cycle lengths lying in the set $S$ an $S$-permutation.  The logarithmic asymptotics of $S$-permutations are governed by the largest element of $S$.

\begin{thm}\label{thm:leading-term} Let $S$ be a finite set of positive integers, with $m = \max S$, and such that the elements of $S$ do not all have a common factor.  Let $n! p_n^{(S)}$ be the number of $S$-permutations of $[n]$.  Then
 \[ p_n^{(S)} n!^{1/m} \sim \cdot C_S n^{-1/2 + 1/2m} \exp(f_S(n^{1/m}))  \]
for some polynomial $f_S$ of degree $m-1$ and constant $C_S$ which can be explicitly computed. In particular, 
\[ \lim_{n \to \infty} {\log p_n^{(S)} \over \log n!} = -1/m. \] \end{thm}
The condition $\gcd S = 1$ is a technical one required so that $\exp (\sum_{s \in S} z^s/s)$ is Hayman-admissible.

\begin{proof} We apply Hayman's method \cite{hayman, wilf} to the generating function $f(z) = \exp \left( \sum_{s \in S} z^s/s \right)$.  We have
 \[ p_n^{(S)} \sim {f(r_n) \over r_n^n \sqrt{2\pi b(r_n)}} \]
where $a(z) = \sum_{s \in S} z^s$, $r_n$ is the positive real root of $a(z) = n$, and $b(z) = \sum_{s \in S} sz^s$.  Using the Lagrange inversion formula, we can find an asymptotic series for $r_n$ in descending powers of $n^{1/m}$.  (See \cite{wilf-bulletin} for details.)  From this we can determine the leading-term asymptotic behavior of $f(r_n)$ and $r_n^n$; we get $f(r_n) = \exp(n/m + c_1 n^{(m-1)/m} + \cdots + c_m n^0 + O(n^{-1/m}))$ and $r_n^n = n^{n/m} \exp(d_1 n^{(m-1)/m} + d_2 n^{(m-2)/m} + \cdots + d_m n^0 + O(n^{-1/m})$ for constants $c_k, d_k$ depending on $S$.  Finally, $b(r_n) \sim mn$.  So
\[ p_n^{(S)} \sim {\exp(n/m + c_1 n^{(m-1)/m} + \cdots + c_m n^0 + O(n^{-1/m}) \over n^{n/m} \exp( d_1 n^{(m-1)/m} + \cdots + d_m n^0 + O(n^{-1/m})} \]
and applying Stirling's approximation gives the result. \end{proof}
To illustrate the theorem, consider $S = \{1, 2, 3\}$, so $r_n$ is the positive real root of $z + z^2 + z^3 = n$.  This has asymptotic series $r_n = n^{1/3} - {1 \over 3} - {2 \over 9} n^{-1/3} + {7 \over 81} n^{-2/3} + O(1/n)$ for large $n$.  From this we can find the leading terms $r_n^n \sim n^{n/3} \exp(-n^{2/3}/3 - 5n^{1/3}/18)$ and $f(r_n) \sim \exp(n/3 + n^{2/3}/6 + 5n^{1/3}/9 - 5/18)$.  Thus
\[ p_n^{(S)} \sim {\exp \left( { n \over 3} + {1 \over 2} n^{2/3} + {5 \over 6} n^{1/3} - {5 \over 18} \right) \over n^{n/3} \sqrt{6\pi n}} \]
and finally
\[ p_n^{(S)} \cdot n!^{1/3} \sim (e^5 2^6 3^9 \pi^6)^{-1/18}  n^{-1/3} \exp \left( {1 \over 2} n^{2/3} + {5 \over 6} n^{1/3} \right) \]

\begin{cor}\label{cor:general-average} The expected number of cycles of length $k$ in an $S$-permutation chosen uniformly at random, where $k \in S$ and $m = \max S$, is $n^{k/m}/k \cdot (1+ o(1))$ as $n \to \infty$. \end{cor}
This has also been shown by Benaych-Georges \cite{benaych-georges} and Timashev \cite{timashev}.
\begin{proof} Let $a_n = n! p_n^{(S)}$ be the number of  $S$-permutations of $[n]$.  The generating function of $S$-permutations by their size and number of $k$-cycles is
\[ G^{(S)}(z,u) = \exp \left( \left( \sum_{s \in S} z^s/s \right) + (u-1) z^k/k \right). \]
The mean number of $k$-cycles in $S$-permutations of $[n]$ is therefore
\[ {[z^n] \left( {\partial \over \partial z} \left. G^{(S)}(z,u) \right|_{u=1} \right) \over [z^n] G^{(S)}(z,1)} = {[z^n] {1 \over k} z^k G^{(S)}(z,1) \over [z^n] G^{(S)}(z,1)} = {1 \over k} {p_{n-k} \over p_n}. \]
Now, $p^{(S)}_{n-1}/p^{(S)}_n \sim (n-1)!^{-1/m}/n!^{-1/m} = n^{1/m}$, the subexponential factor in Theorem \ref{thm:leading-term} being slowly varying.  So the mean number of $k$-cycles is asymptotic to ${1 \over k} (n^{1/m})^k$, as desired.
\end{proof}

The Boltzmann sampler \cite{dfls} for $S$-permutations provides an explanation for Corollary \ref{cor:general-average}.  To generate random $S$-permutations, we fix a positive real parameter $x$ and then pick a cycle type by taking $\Po(x^k/k)$ cycles of length $k$ for each $k \in S$. The cycles themselves are then populated with elements uniformly at random.  Fixing $x$ to be the positive root of $\sum_{k \in S} z^{k} = n$ -- that is, $x = r_n$ -- gives permutations of expected size $n$, and all $S$-permutations of the same size are equally likely to be generated.  The expected number of $k$-cycles of a permutation generated by this process is $r_n^k/k$.  

Alternatively, we could find the number of permutations with all cycle lengths in some set $S$ with largest element $m$ by summing over cycle types.  For example, for involutions we have the sum $\sum_{l+2k = n} {n! \over l! k! 2^k}$, or $\sum_{k=0}^{n/2} {n! \over (n-2k)! k! 2^k}$.   More generally, if $S = \{ s_1, \ldots, s_j \}$ with $s_1 < s_2 < \ldots < s_j$, then the number of $S$-permutations of $[n]$ is given by
\begin{equation}\label{eq:cycle-type-sum} \sum_{c_1 s_1 + \ldots + c_j s_j = n} {n! \over c_1! \ldots c_j! s_1^{c_1} \ldots s_j^{c_j}} \end{equation}
where the sum is over $j$-tuples of positive integers $(c_1, \ldots, c_j)$.  This summand can be approximated as a multivariate Gaussian integral.  It may be possible to do this integral using Laplace's method; Greenhill et al. \cite[Thm. 6.4]{greenhill-janson-rucinski} give a version of this method adapted to sums over high-dimensional lattices.  We do not carry out this computation as the method of Theorem \ref{thm:leading-term} is effective.

  \subsection{Graph-theoretic decomposition}
An involution $\sigma$ can be represented as a partial matching on the set $[n]$, where $k$ and $l$ are matched if $\sigma(k) = l$ (and therefore $\sigma(l) = k$).   We can view this matching as a graph, by drawing an edge between $k$ and $l$ when $\sigma(k) = l$.   A pair of partial matchings or involutions, $(\sigma, \tau)$, can be identified with a graph on the vertex set $[n]$ with 2-colored edges, where we color the edges solid or dotted according to whether they are from $\sigma$ or from $\tau$.  We write $\sigma \cup \tau$ for this graph, and refer to it as as a {\it superposition}, and $\tau \circ \sigma$ for the corresponding permutation.

\begin{thm}\label{thm:trivariate}
The trivariate generating function for pairs of partial matchings $(\sigma, \tau)$, counted according to the size of the ground set (indicated by the variable $z$) and number of paths and cycles  in $\sigma \cup \tau$ (indicated by $u$ and $v$ respectively), exponential in $z$ and ordinary in $u$ and $v$, is
\begin{equation}\label{eq:involution-pair-trivariate} Q(z,u,v) = {\exp(uz/(1-z)) \over (1-z^2)^{v/2}} \end{equation}
That is, $n! [z^n u^k v^l] Q(z,u,v)$ is the number of pairs of partial matchings on $n$ vertices with $k$ paths and $l$ cycles.
\end{thm}

\begin{proof}
We enumerate the possible connected components of a pair of partial matchings and apply the exponential formula.

The connected components of such a graph are cycles of even length and paths, with the edges alternating in color.  These are the only possible components since if colors are ignored, all vertices must have degree at most two.  We note that the degenerate path (a single vertex) and the degenerate cycle (two vertices connected by a solid edge and a dotted edge) are both possible components.  We consider the length of a path to be its number of vertices, so a single vertex is a path of length 1. 

For $n \ge 3$, the number of labelled cycles of length $n$ is $((n-1)!)/2$; there are two ways to color a labelled cycle of even length with alternating edge colors, and zero ways to color a cycle of odd length.  There is exactly one labelled cycle on two vertices.  Thus the exponential generating function (egf) for properly colored cycles is 
\[ 1! {z^2 \over 2!} + 3! {z^4 \over 4!} + 5! {z^5 \over 6!} + \cdots = {z^2 \over 2} + {z^4 \over 4} + {z^6 \over 6} + \cdots = {1 \over 2} \log {1 \over 1-z^2}. \]
The number of labelled paths of length $n$ is $n!/2$ for $n \ge 2$, and $1$ for $n = 1$.  There are two ways to color a path of any length -- pick a color for a single edge and the coloring can be completed in exactly one way -- except that there is only one way to color a 1-path.  Thus there are $n!$ colored, labelled $n$-paths, for each $n \ge 1$, and the egf for properly colored paths is $z/(1-z)$.

The generating function of components marked according to their type (path or cycle) is therefore $u \cdot {z \over 1-z} + v \cdot {1 \over 2} \log {1 \over 1-z^2} $, and applying the exponential formula gives (\ref{eq:involution-pair-trivariate}).
\end{proof}

We quickly derive two corollaries more relevant to permutation enumeration.

\begin{cor} The exponential generating function of pairs of involutions is $P(z) = \exp(z/(1-z))/\sqrt{1-z^2}$. \end{cor}
\begin{proof}  Take the specialization $u = 1, v = 1$ in Theorem \ref{thm:trivariate}.  This gives the exponential generating function for pairs of partial matchings, which we identify with pairs of involutions. \end{proof}

\begin{figure}
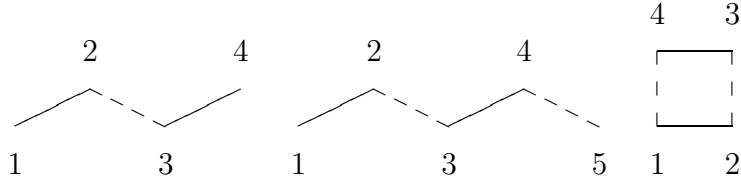
\label{fig:paths}
\[ 
\xy
(10,5)*{1};
(10,10)*{}; (20,15)*{} **\dir{-};
(20,20)*{2};
(20,15)*{}; (30,10)*{} **\dir{--};
(30,5)*{3};
(30,10)*{}; (40,15)*{} **\dir{-};
(40,20)*{4};
\endxy
\: \: \: \: \: \:
\xy
(10,5)*{1};
(10,10)*{}; (20,15)*{} **\dir{-};
(20,20)*{2};
(20,15)*{}; (30,10)*{} **\dir{--};
(30,5)*{3};
(30,10)*{}; (40,15)*{} **\dir{-};
(40,20)*{4};
(40,15)*{}; (50,10)*{} **\dir{--};
(50,5)*{5};
\endxy
\: \: \: \: \: \:
\xy
(10,5)*{1};
(10,10)*{}; (20,10)*{} **\dir{-};
(20,5)*{2};
(20,10)*{}; (20,20)*{} **\dir{--};
(10,25)*{4};
.(20,20)*{}; (10,20)*{} **\dir{-};
(20,25)*{3};
(10,20)*{}; (10,10)*{} **\dir{--};
\endxy
\]
\caption{Paths of length four and five, and a cycle of length four.}
\end{figure}

\begin{cor}\label{cor:inv-bivariate} The semi-exponential generating function of pairs of involutions $(\sigma, \tau)$, counted by the size of the ground set and the number of {\it permutation} cycles in the composition $\tau \circ \sigma$, is \[ R(z,u) = {\exp \left( {uz \over 1-z} \right) \over (1-z^2)^{u^2/2}}, \] \end{cor}
\begin{proof}  Consider a pair of perfect matchings $(\sigma, \tau)$.  Each connected component of the corresponding graph $\sigma \cup \tau$ gives rise to either one or two cycles in $\tau \circ \sigma$.  Each $2k$-cycle in the graph $\sigma \cup \tau$ gives rise to two $k$-cycles in the permutation $\tau \circ \sigma$, corresponding to half of the vertices.  Each $k$-path in $\sigma \cup \tau$ gives rise to a $k$-cycle in $\tau \circ \sigma$.  For example, the paths illustrated in Figure \ref{fig:paths} correspond to the permutation cycles $(1342)$ and $(13542)$ respectively, and the cycle to the pair of permutation cycles $(13)(24)$.  To count by permutation cycles, then, we need $z^n u^k v^l$ in $Q(z,u,v)$ to be mapped to $z^n u^{k+2l}$; thus we take the specialization $R(z,u) = Q(z,u,u^2)$ in Theorem \ref{thm:trivariate}.  \end{proof}

  \subsection{Asymptotic distribution of the number of k-cycles}

In this section we show 

\begin{thm}\label{thm:k-cycle-dist}The distribution of the number of $k$-cycles of the composition of a pair of random involutions of $[n]$ converges in distribution to the distribution of $A_k + 2B_k$ as $n \to \infty$, where $A_k$ and $B_k$ are independent, $A_k$ is Poisson of mean 1, and $B_k$ is Poisson of mean $1/(2k)$.\end{thm}

We need the following more general result.  Recall that a sequence $\{ b_n \}$ is {\it slowly varying} if $\lim_{n \to \infty} b_{n-1}/b_n = 1$.
\begin{lem}\label{lem:factorial-mgf}
Let $P(z,u)$ be the generating function 
\[ P(z,u) = \sum_{n,k \ge 0} P_{n,k} {z^n \over n!} u^k \]
where $P_{n,k}$ is the number of objects in a combinatorial class $\mathcal{P}$ of size $n$ with some parameter $\chi$ equal to $k$.  Assume $P(z,u) = Q(z) e^{R(z,u)}$ with $R$ a polynomial, $[z^n] Q(z)$ is slowly varying as $n \to \infty$, and $R(1,t)$ is the factorial moment generating function of some distribution which is determined by its moments.  For each $n$, let $\P_n(\chi=k) = P_{n,k}/\sum_k P_{n,k}$ define a probability distribution on the positive integers.  Then as $n \to \infty$, the sequence of distributions of $\chi$ on $\P_n$ converges in distribution to the distribution with factorial moment generating function $\exp R(1,t)$.
\end{lem}
\begin{proof}
Let $j$ be the degree of $R$ in the variable $u$.  We can show by induction that $\partial_u^r P(z,1) = P(z,1) T_r(z)$, where $T_r(z)$ has degree $jr$.  Then we have
\[ \E_n((\chi)_r) = {[z^n] \left. \partial_u^r P(z,u) \right|_{u=1} \over [z^n] P(z,1)} = {[z^n] P(z,1) T_r(z) \over [z^n] P(z,1)}. \]
Now, $\lim_{n \to \infty} {[z^{n-s}] P(z,1)/[z^n] P(z,1)} = 1$ from the condition on slow variation.  Then
\[ {[z^n] P(z,1) T_r(z) \over [z^n] P(z,1)} = \sum_{s=0}^{jr} [z^s] T_r(z) {[z^{n-s}] P(z,1) \over [z^n] P(z,1)} \]
and taking limits as $n \to \infty$ gives 
\[ \lim_{n \to \infty} {[z^n] P(z,1) T_r(z) \over [z^n] P(z,1)} = T_r(1). \]
So $\lim_{n \to \infty} \E_n((\chi)_r) = T_r(1)$.  Now let $F(t) = R(1,t)$.  The $r$th factorial moment of the distribution with factorial mgf $F(t)$ is $F^{(r)}(1)$.  This is $ \left. {\partial^r \over \partial u^r} \exp R(z,u) \right|_{z=u=1} $
and we recall that $\partial_u^r Q(z) \exp R(z,u)|_{u=1} = P(z,1) T_r(z)$ by definition.  Therefore we have
\[ \left. {\partial^r \over \partial u^r} \exp R(z,u) \right|_{z=u=1} = {P(1,1) T_r(1) \over Q(1)} = T_r(1) \]
which is what we wanted. \end{proof}

\begin{proof}[Proof of Theorem \ref{thm:k-cycle-dist}.] We apply Lemma \ref{lem:factorial-mgf} with 
 \[ P(z,u) = {\exp(z/(1-z)) \over \sqrt{1-z^2}} \exp \left( (u-1) z^k + (u^2-1) {z^{2k} \over 2k} \right). \]
The slow variation hypothesis holds since $[z^n] \exp(z/(1-z))/\sqrt{1-z^2} = a_n^2/n! = e^{2\sqrt{n}}(8\pi en)^{-1/2}(1+o(1))$.  We have $R(z,u) = \exp ((u-1)z^k + (u^2-1) z^{2k}/(2k))$; this is the factorial moment generating function of $A_k + 2B_k$, which follows from the fact that $Poisson(\lambda)$ has factorial mgf $\exp(\lambda(t-1))$.  Finally, we recall that if the moment generating function of a random variable has positive radius of convergence, then the random variable is determined by its moments \cite[Thm 30.1]{billingsley}.  $A_k + 2B_k$ has mgf $\exp(e^t-1 + (e^{2t}-1)/2k)$, which is entire.
\end{proof}

The sum of Poissons given in Theorem \ref{thm:k-cycle-dist} is quite natural.  There are two types of components in superpositions of partial matchings on $[n]$ that can lead to $k$-cycles of the corresponding permutations: paths of length $k$ (which induce one permutation $k$-cycle) and cycles of length $2k$ (which induce two permutation $k$-cycles).  For large $n$ and fixed $k$, the expected number of $k$-paths approaches 1 and the expected number of $2k$-cycles approaches $1/k$.  Furthermore, the sites in which individual cycles can appear are each rare, so it is not surprising to see an independent Poisson distribution for each type of component.

\subsection{Partial matchings with a specified number of fixed points}

In this section we consider superpositions of partial matchings, $\sigma \cup \tau$, where $\sigma$ is chosen uniformly from all partial matchings on $[n]$ with $k$ fixed points, and $\tau$ is chosen uniformly from all partial matchings with $l$ fixed points.

\begin{prop} \label{prop:exp-r-paths} The expected number of $r$-paths in $\sigma \cup \tau$ is
 \begin{equation}\label{eq:odd-r-paths} kl {({n-k \over 2})_{(r-1)/2} ({n-l \over 2})_{(r-1)/2} 2^{r-1} \over (n)_r} \end{equation}
if $r$ is odd, and
\begin{equation}\label{eq:even-r-paths} (k(k-1)({n-k \over 2})_{r/2-1}({n-l \over 2})_{r/2} + l(l-1)({n-k \over 2})_{r/2} ({n-l \over 2})_{r/2-1}) 2^{r-1} \over 2(n)_r \end{equation}
if $r$ is even. \end{prop} 

\begin{prop}\label{prop:exp-r-cycles} The expected number of $r$-cycles in $\sigma \cup \tau$, is
\begin{equation}\label{eq:r-cycles} { \left( {n-k \over 2} \right)_{r/2} \left( {n-l \over 2} \right)_{r/2} 2^r \over r(n)_r} \end{equation}
for even $r$. \end{prop}

These statements can be easily verified.  For odd paths, we compute the probability that a path occurs which traverses the edges $1, 2, \ldots, r$ in that order, and multiply by the number of possible paths.  The argument is similar for even paths, except we must handle the cases where the two ends of the path are fixed points in $\sigma$ and fixed points in $\tau$ separately.  Finally, we do this for cycles; the most interesting feature is the factor of $r$ in the denominator which arises from the symmetry of cycles.  This model of random involutions, with $n \to \infty$ and $k, l$ varying with $n$ in such a way that $k+l = \Omega(1)$ and $k+l = o(n)$ simultaneously, has been considered in \cite{roberts-vivaldi} in the context of dynamical systems.

\begin{cor}The expected number of paths of length $r$ in $\sigma \cup \tau$, the superposition of two randomly selected perfect matchings on $[n]$, where $\sigma$ and $\tau$ each have $pn$ fixed points, is asymptotic to $p^2 (1-p)^{r-1}$ as $n \to \infty$.  The expected number of cycles of length $r$ (if $r$ is even) approaches $(1-p)^r/r$ as $n \to \infty$. \end{cor}

\begin{cor}\label{cor:sqrt-cycles} Let $r = O(\sqrt{n})$  as $n \to \infty$.  Then the mean number of $r$-paths in $\sigma \cup \tau$, where $\sigma$ and $\tau$ are randomly selected involutions with $\sqrt{n}$ fixed points each, is asymptotic to $\exp(-r/\sqrt{n})$ as $n \to \infty$, and the mean number of $r$-cycles is asymptotic to $\exp(-r/\sqrt{n})/r$. \end{cor}

These follow from Propositions \ref{prop:exp-r-paths} and \ref{prop:exp-r-cycles} by making appropriate substitutions, and applying Stirling's formula in the case of Corollary \ref{cor:sqrt-cycles}.  In both cases the number of $r$-paths decays exponentially in $r$.  This can be explained by considering the process of path formation as a random walk from one fixed point to another on the complete graph $K_n$.  The length of a path is the time such a walk takes to reach some fixed point, which is geometric and has expectation the inverse of the proportion of points which are fixed points.

The last of these results can be translated back into the terminology of involutions.  To get a better sense of the scaling behavior of cycle sizes, we look at the expected number of $\alpha\sqrt{n}$-cycles of a composition of two random involutions.  The expected number of $k$-cycles is $\left(b_{n-k} + {1\over 2k} b_{n-2k}\right)/b_n$, where $b_n = [z^n] \exp(z/(1-z))/\sqrt{1-z^2}$.  Recall that $b_n \sim e^{2\sqrt{n}} (8\pi e n)^{-1/2}$.  Let $k = \alpha \sqrt{n}$ grow with $n$.  One can compute that $b_{n-\alpha \sqrt{n}} \sim e^{-\alpha}$ as $n \to \infty$.  Such square-root scaling is typical of structures counted by generating functions which are the exponential of a function with a simple pole.  The simplest example may be the ``fragmented permutations'' (permutations with rooted cycles) or ``sets of lists'' counted by $\exp(z/(1-z))$, which will be treated in more detail in \cite{lugo-further}.

\subsection{The total number of cycles}

The function $R(z,u)$ given in Corollary \ref{cor:inv-bivariate} will be our jumping-off point for asymptotic results on cycle structure.  We need the following asymptotic result.

\begin{thm}[Wright]\cite[Thm. 2 and Thm. 3]{wright-32}\label{thm:wright}
\begin{enumerate} \renewcommand{\labelenumi}{(\alph{enumi})}
 \item  The leading-term asymptotics for
\[ c_n = [z^n] (1-z)^{\beta} \Phi(z) \exp \left( {1 \over 1-z} \right) \]
where $\beta$ is a complex number and $\Phi(z)$ is regular in the unit disk are given by  
\[ c_n = {1 \over n^{\beta/2 + 3/4}} \left[ \exp(2\sqrt{n}) {1 \over 2 \sqrt{\pi}} \Phi(1) e^{1/2} \right] (1 + O(n^{-1/2})). \]
\item The leading-term asymptotics for 
\[ [z^n] \left( \log{1 \over 1-z}\right)^k (1-z)^{\beta} \Phi(z) \exp \left( {1 \over 1-z} \right) \]
with $k$ a positive integer can be derived from that for the $k=0$ case by differentiating $k$ times with respect to $\beta$ and switching signs if $k$ is odd. \end{enumerate} 
\end{thm}

In particular, applying (a) with $\Phi(z) = e^{-z}/\sqrt{1+z}, \beta = -1/2$ gives $[z^n] P(z) = {1 \over \sqrt{8\pi e n}} \exp(2\sqrt{n})$.  This is consistent with the known number of involutions in (\ref{eq:number-of-involutions}).

\begin{prop}\label{prop:mean-number-paths} The mean number of components which are paths in a superposition of two partial matchings on $[n]$ selected uniformly at random is $\sqrt{n} + O(1)$. \end{prop}
\begin{proof} The bivariate generating function counting superpositions of partial matchings by size and number of paths is $Q(z,u,1) = \exp(uz/(1-z))/\sqrt{1-z^2}$.  A standard result in generating functions is that the expectation of the parameter marked by $u$ in such a bivariate generating function is $[z^n] Q_u(z,u,1)/[z^n] Q(z,u,1)$ for objects of size $n$.  Differentiating and applying Theorem \ref{thm:wright} gives the result. \end{proof}

\begin{prop}\label{prop:mean-number-cycles} The mean number of components which are cycles in a superposition of two partial matchings on $[n]$ selected uniformly at random is ${1 \over 4} \log{n} + O(n^{-1/2} \log n)$. \end{prop}
\begin{proof} The bivariate generating function counting superpositions of partial matchings by size and number of cycles is $Q(z,1,v) = \exp(z/(1-z)) (1-z^2)^{-v/2}$; differentiate and apply Theorem \ref{thm:wright}.
\end{proof}

\begin{prop} The mean number of elements in cycles in a superposition of two random partial matchings of $[n]$ is ${1 \over 2} \sqrt{n} + O(1)$. \end{prop}
\begin{proof} The generating function counting pairs of matchings by their size and number of elements in cycles is
 \[ \exp \left( {z \over 1-z} + {u^2 z^2 \over 2} + {u^4 z^4 \over 4} + {u^6 z^6 \over 6} + \cdots \right) = {\exp(z/(1-z)) \over \sqrt{1-u^2 z^2}}. \]
Again, we differentiate and apply Theorem \ref{thm:wright}.
\end{proof}

\begin{thm}\label{thm:mean-number-permutation-cycles} The mean number of cycles in a composition of two uniform random involutions on $[n]$ is $\sqrt{n} + {1 \over 2} \log n + O(1)$. \end{thm}
\begin{proof} A superposition of partial matchings with $k$ paths and $l$ (graph) cycles is identified with a composition of involutions having $k + 2l$ (permutation) cycles.   \end{proof}

\begin{prop} The probability that a superposition of two partial matchings of $[n]$ selected uniformly at random has no cyclic components is $\sqrt{2} n^{-1/4} + O(n^{-3/4})$ as $n \to \infty$. \end{prop}
\begin{proof} Partial matchings with no cyclic components have generating function $Q(z,1,0) = \exp(z/(1-z))$; thus the probability in question is
 \[ {[z^n] \exp \left( {z \over 1-z} \right) \over [z^n] \exp \left( {z \over 1-z} \right)/\sqrt{1-z^2}} \]
By Theorem \ref{thm:wright} the numerator is $e^{2\sqrt{n}}/(2n^{3/4} \sqrt{e\pi}) (1 + O(n^{-1/2}))$; the denominator is $e^{2\sqrt{n}}/\sqrt{8\pi e n} (1+O(n^{-1/2}))$, giving the desired result. \end{proof}

\subsection{Fixed-point-free involutions}

\begin{prop}\label{prop:a} The number of pairs of fixed-point-free involutions $(\sigma, \tau)$ on $[2n]$ such that $\pi = \tau \circ \sigma$ has $2c_k$ $k$-cycles for each $k$ is the same as the number of permutations of $[2n]$ which have $c_k$ $2k$-cycles for each $k$, and no cycles of odd length. \end{prop}

\begin{proof} We construct a bijection between the two sets.  Given such a pair of fixed-point-free involutions, the graph of $\sigma \cup \tau$ consists of $c_k$ graph cycles of length $2k$, with the edges alternately solid and dotted.  From each graph cycle we construct a permutation cycle.  We need only make a choice of direction, say by starting at the smallest element and following the solid edge out of that element.  This operation is clearly reversible; given a permutation with only even cycles we can reconstruct the graph $\sigma \cup \tau$ of a pair of fixed-point-free involutions. \end{proof}

\begin{prop} The number of cycles in a composition of two fixed-point-free involutions on $[2n]$ chosen uniformly at random has the distribution of $2\sum_{k=1}^n X_k$, where $X_k$ is Bernoulli with mean $1/(2k-1)$ and the $X_k$ are independent. \end{prop}

\begin{proof} The distribution of the number of cycles of a permutation of $[2n]$ with all cycle lengths even is that of $\sum_{k=1}^n X_k$, where $X_k$ is Bernoulli with mean $1/(2k-1)$ and the $X_k$ are independent \cite[Thm. 3.7]{lugo}.  From Proposition \ref{prop:a}, there are exactly the same number of permutations of $[2n]$ with $2j$ cycles, all of even length, as there are pairs of fixed-point-free involutions $(\sigma, \tau)  \in S_{2n} \times S_{2n}$ with $\tau \circ \sigma$ having $j$ cycles. \end{proof}

Note that the expected number of cycles in a composition of two fixed-point-free involutions of $[2n]$ is $2H_{2n} - H_n = \log n + (2 \log 2 + \gamma) + O(n^{-2})$, which differs from the expected number of cycles in a random permutation of $[2n]$ by $\log 2 + O(n^{-1})$.   However, compositions of fixed-point-free involutions do not ``look like'' random permutations.  Most obviously, a composition of fixed-point-free involutions of $[n]$ has no cycles longer than $n/2$. Cycle lengths satisfy the following limit law.

\begin{prop} Fix constants $0 \le \gamma \le \delta \le 1/2$.  Let $p_i(n;\gamma,\delta)$ be the probability that $1$ is contained in a cycle of $\tau \circ \sigma$ of length between $\gamma n$ and $\delta n$, where $\sigma$ and $\tau$ are fixed-point-free involutions on $[n]$ chosen uniformly at random. Then $\lim_{n \to \infty} p_i(n;\gamma,\delta) = \sqrt{1-2\gamma}-\sqrt{1-2\delta}$. \end{prop}
\begin{proof} Call a permutation with all cycle lengths even an $E$-permutation, and a composition of fixed-point-free involutions an $I$-permutation.  The number of $2k$-cycles in $E$-permutations of $[n]$ is half the number of $k$-cycles in $I$-permutations of $[n]$, by Proposition \ref{prop:a}.  In particular the number of elements of $2k$-cycles in $E$-permutations of $[n]$ and the number of elements of $k$-cycles in $I$-permutations of $[n]$ are equal.  So the probability that a random element of a random $I$-permutation of $[n]$ is in a cycle of length in $[\gamma n, \delta n]$ is equal to the probability that a random element of a random $E$-permutation of $[n]$ is in a cycle of length in $[2\gamma n,  2 \delta n]$.  By \cite[Thm. 3.5]{lugo} the latter probability approaches $\sqrt{1-2\gamma} - \sqrt{1-2\delta}$ as $n \to \infty$.  \end{proof}

\begin{prop} Fix $\epsilon \in (0,1/2)$.  The expected number of elements in $k$-cycles in a composition of two random fixed-point-free involutions of $[n]$ converges uniformly to $(1-2k/n)^{-1/2}$ as $k/n \to \infty$ with $0 < k/n < 1/2-\epsilon$.\end{prop}
\begin{proof}By Proposition \ref{prop:exp-r-cycles} the expected number of elements in $r$-cycles in a superposition $\sigma \cup \tau$ of fixed-point-free perfect matchings is
\[ {(n/2)!^2\over((n-r)/2)!^2} 2^r {(n-r)! \over n!}. \]
In the case $r = \alpha n$, this is asymptotic to $1/\sqrt{1-\alpha}$ as $n \to \infty$, with uniform convergence over $0 < \alpha < 1$; this is shown in \cite[Prop. 3.4]{lugo}, where the same expression occurs in relation to permutations with all cycle lengths even.  Noting that elements in $r$-cycles in a pair of perfect matchings give rise to elements in $r/2$-cycles of the corresponding permutation gives the desired result. \end{proof}

  \subsection{The number of factorizations of a permutation into involutions}

The square of the number of involutions of $[n]$ is a bit larger than $n!$; we have
\[ a_n^2 \sim n! \cdot {e^{2 \sqrt{n}} \over \sqrt{8\pi e n}}. \]
The mean number of factorizations of a random permutation into a product of involutions is just the second factor.  The number of factorizations can be as large as $a_n$ for the identity permutation, since $id = \sigma^2$ for any involution $\sigma$, and as small as $n-1$ for those permutations which consists of an $(n-1)$-cycle and a $1$-cycle.

\begin{thm}\label{thm:inv-factorizations} Define the function 
\[ f(r,k) = \sum_{j=0}^{\lfloor r/2 \rfloor} {r! \over (r-2j)! j! 2^j} k^{r-j}. \]
Let $\pi$ be a permutation of $[n]$ with $c_k$ cycles of length $k$, for each $k$.  Then 
\[ F(\pi) = \prod_{k=1}^n f(c_k, k) \] is the number of factorizations of $\pi$ into two involutions, i. e. the number of solutions of $\pi = \tau \circ \sigma$ with $\sigma$ and $\tau$ involutions. \end{thm}

We remark that $f(r,k)$ is the number of partial matchings of $[r]$ with $k$-colored components.  This interpretation is key to the proof, which works by pairing up some of the $k$-cycles and then assigning one of $k$ partial factorizations to each unpaired $k$-cycle or pair of $k$-cycles.

We begin with the following special case.

\begin{lem}\label{lem:graphical-path} The number of ways to factor an $n$-cycle $\pi$ into two involutions is $n$. \end{lem} 
\begin{proof} Without loss of generality let $\pi = (123\cdots n)$.  We construct a corresponding pair of partial matchings $(\sigma, \tau)$.  This must be a path of length $n$, since cycles in $\sigma \cup \tau$ give rise to pairs of permutation cycles.  So we consider an unlabeled path of length $n$ with alternating solid and dotted edges, and attempt to label it.  We begin by labelling some vertex by $1$.  Then follow the solid edge at that vertex, followed by the dotted edge at the next vertex, to determine the site of 2; repeat to determine the sites of 3, 4, and so on.  The remaining vertices can therefore be labelled in exactly one way.  \end{proof}

For example, the cycle $(1234)$ has the factorizations 
\[ (\sigma, \tau) = ((1)(24)(3), (12)(34)), ((13)(2)(4), (14)(23)), ((12)(34), (2)(13)(4)), ((14)(23), (24)(1)(3)). \] 

This is also a special case of a formula given in \cite[Thm. 2.1]{goupil-schaeffer} for the number of factorizations of an $n$-cycle into permutations of types $\lambda$ and $\mu$.  Note that if $n = 2k+1$, then $\lambda$ and $\mu$ each have type $2^k 1$; if $n = 2k$, either $\lambda$ has type $2^{k-1} 1^2$ and $\mu$ has type $2^k$, or vice versa.

\begin{lem}\label{lem:graphical-cycle} The number of ways to factor a permutation $\pi$ of $[2n]$ consisting of two $n$-cycles into two involutions $\sigma, \tau$, such that the corresponding graph $\sigma \cup \tau$ is a $2n$-cycle, is $n$. \end{lem}
\begin{proof} 
Without loss of generality, let $\pi = (1, 2, \ldots, n) (n+1, n+2, \ldots, 2n)$ in cycle notation.  We draw a graphical cycle with $2n$ vertices, with edges alternately solid and dotted.  Label some arbitrary vertex with $1$; follow solid and dotted edges alternately around the cycle to place $2, 3, \ldots, n$.  Then label some arbitrary unlabeled vertex with $n+1$ and follow solid and dotted edges alternately around the cycle to place $n+2, \ldots, 2n$.  There are $2n^2$ ways to carry out this procedure.  However, the unlabeled $2n$-cycle with alternately colored edges has $2n$ symmetries.  So there are $(2n^2)/(2n) = n$ distinct labellings; each one corresponds to a factorization. 
\end{proof}

\begin{proof}[Proof of Theorem \ref{thm:inv-factorizations}.]  Given an arbitrary permutation $\pi$ to be factored into involutions with $\pi = \tau \circ \sigma$, we can consider the cycles of each length separately.  Consider the cycles of length $k$; assume there are $r$ of these.  We pair up some of the $k$-cycles with each other, representing that they come from the same cycle in the graph $\sigma \cup \tau$.  Those cycles which remain unpaired arise from paths, not cycles, in $\sigma \cup \tau$.  We then factor each unpaired cycle according to  Lemma \ref{lem:graphical-path}, and each pair of cycles according to Lemma \ref{lem:graphical-cycle}.  If there are $j$ pairs of cycles, then there are $r-2j$ unpaired cycles, and thus $r-j$ total components to factor; thus the number of such factorizations, once the cycles are paired up, is $k^{r-j}$.  The number of ways to find $j$ disjoint pairs of cycles, with order irrelevant, is
\[ {{r \choose 2} {r-2 \choose 2} \ldots {r - 2j+2 \choose 2} \over j!} = {r! \over (r-2j)! j! 2^j}. \]  Summing over $j$ gives the function $f(r,k)$ defined in the theorem. \end{proof}

It appears that the distribution of the number of factorizations of a random permutation of $[n]$ into involutions approaches a lognormal distribution as $n \to \infty$.  We consider the following probability model: let $X_k = \Po(1/k)$ for $k = 1, 2, \ldots, n$, where $n$ is a positive integer parameter.  Let $m = X_1 + 2X_2 + \cdots + nX_n$, and take a permutation of $[m]$ with $X_k$ cycles of length $k$ for each $k$, chosen uniformly at random from all permutations of that cycle type.  We denote the corresponding measure on the set of all permutations by $\P^*_n$. Then $\E \left( \sum_{k=1}^n kX_k \right) = n$.  This model generates each permutation of $[m]$ having all cycle lengths less than or equal to $n$ with probability $e^{-H_n}/m!$, where $H_n = \sum_{k=1}^n 1/k$ is a harmonic number; in particular for each $m \le n$, each permutation of $[m]$ occurs with the same probability.  

\begin{thm}\label{thm:lognormal} As $n \to \infty$,
\[ \lim_{n \to \infty} \P^*_n \left( {\log(F(\pi)) - {1 \over 2} (\log n)^2 \over {1 \over 3} (\log n)^3} \le x \right) \to \Phi(x) \] \end{thm}
\begin{proof} First, we show that $\E(\log f(X_k, k)) = \log(k)/k + O(k^{-3})$.  Let $\mu_k = \E(\log f(X_k,k))$.  We can write the expectation as a sum over possible values of $X_k$, giving
\begin{equation}\label{eq:mu-k-sum} \mu_k = e^{-1/k} \sum_{r \ge 1} {1 \over r! k^r} \log f(r,k). \end{equation}
We can derive an asymptotic series for $\log f(r,k)$ from the Taylor series for $\log(1+x)$ around $x = 0$; this gives an asymptotic series for the $r$th term in (\ref{eq:mu-k-sum}), which is of order $k^{-r} \log k$.  Adding these gives
$\mu_k = (\log k)/k + (1/2k^3) + O(k^{-4})$.  Similarly let $h_k = \E((\log f(X_k,k))^2)$; then in like manner we can derive the series
$h(k) = {(\log k)^2/k} + {(\log k)^2/k^2} + {2 \log k/k^3} + O(k^{-4}).$  The variance is given by $\sigma_k^2 = \V(\log f(X_k,k)) = h(k) - \mu_k^2$, and we find $\sigma_k^2 = {(\log k)^2/k} + {2 \log k/k^3} + O(\log k / k^4)$.

Next we show that $\sum_{k=1}^n \mu_k \sim (\log n)^2/2$ and $\sum_{k=1}^n \sigma_k^2 \sim (\log n)^3/3$ as $n \to \infty$.  We have $\mu_k = (\log k)/k + O(k^{-3})$.  Now, $\sum_{k=1}^n (\log k)/k \sim \int_1^n (\log k)/k \: dk = {1 \over 2} (\log n)^2$, where the asymptotic equality can be justified by the Euler-Maclaurin summation formula.  Expanding the big-$O$ notation, $|\mu_k - (\log k)/k| \le Ck^{-3}$ for some constant $C$, so $\sum_{k=1}^\infty \mu_k - (\log k)/k$ converges.  Therefore $\sum_{k=1}^n \mu_k \sim \sum_{k=1}^n (\log k)/k \sim {1 \over 2} (\log n)^2$.  The proof for $\sum_{k=1}^n \sigma_k^2$ is similar.  Note that $\sum_{k=1}^n \mu_k = \E (\log F(\pi))$ and $\sum_{k=1}^n \sigma_k^2 = \V(\log F(\pi))$.
Finally, we apply Lyapunov's central limit theorem to show that $\log F(\pi)$ is asymptotically normal.  We recall the theorem: let $Y_1, Y_2, \ldots$ be independent random variables with finite mean and variance, $\E(Y_n) = \mu_n$ and $\V(Y_n) = \sigma_n^2$.  Let $s_n^2 = \sum_{k=1}^n \sigma_k^2$.  If for some $\delta > 0$, $\E(|Y_k|^{2+\delta})$ is finite for $k = 1, 2, \ldots$ and the Lyapunov condition
\begin{equation}\label{eq:lyapunov-condition} \lim_{n \to \infty} {1 \over s_n^{2+\delta}} \sum_{k=1}^n \E(|Y_k - \E Y_k|^{2+\delta}) = 0 \end{equation}
is satisfied, then the standardization $(\sum_{k=1}^n (Y_n-\mu_n))/s_n$ converges in distribution to a standard normal random variable as $n \to \infty$. We will take $\delta = 1$, and $Y_k = \log f(X_k,k)$.  As previously shown, $s_n^2 \sim (\log n)^3/3$, so $s_n^3 \sim (\log n)^{9/2}/(3\sqrt{3})$.  We also observe $\E(|Y_k|^3)$ is finite for each $k$.  To check (\ref{eq:lyapunov-condition}), first note that
\[ \E(|Y_k - \E Y_k|^3) = \sum_{r \ge 1} \left[ ( \log f(r,k) - \E Y_k)^3 \P(X_k = r) \right] + \left( \E Y_k \right) \P(X_k = 0). \]
Since $\E Y_k$ is positive, this is less than
\[ \left[ \sum_{r \ge 1} (\log f(r,k))^3 \P(X_k = r) \right] + \left( \E Y_k \right) \P (X_k = 0). \]
The first term in this equation is in fact $\E(Y_k^3)$.  (The sum giving $\E(Y_k^3)$ should naturally be over $r \ge 0$, but $f(0,k) = 1$ and so the $r = 0$ term does not contribute to the sum.)  Therefore we have
\[ \E(|Y_k - \E Y_k|^3) \le \E(Y_k^3) + (\E Y_k) \P(X_k = 0) \le \E(Y_k^3) + \E(Y_k). \]
But $\E(Y_k^3) \sim (\log k)^3/k$ and $\E Y_k \sim (\log k)/k$ as $k \to \infty$.  so $\E(|Y_k - \E Y_k|^3) \sim (\log k)^3/k$.  Therefore we have
\[ {1 \over s_n^3} \sum_{k=1}^n \E(|Y_k - \E Y_k|^3) \sim {3^{3/2} \over (\log n)^{9/2}} {(\log n)^4 \over 4} = {3^{3/2}/4 \over \sqrt{\log n}} \]
and in particular this goes to $0$ as $n \to \infty$, so (\ref{eq:lyapunov-condition}) is satisfied.  Therefore the standardization of $\log F(\pi)$ converges in distribution to the standard normal, as desired. \end{proof}

Simulation experiments lead to the following conjecture.

\begin{conj} Let $\pi$ be a permutation of $[n]$ chosen uniformly at random.  There exists a positive constant $c \approx 0.16$, such that
\[ \lim_{n \to \infty} \P \left( {\log(F(\pi)) - {1 \over 2} (\log n)^2 \over c (\log n)^3} \le x \right) \to \Phi(x) \]
where $\Phi(x)$ is the distribution function of the standard normal. \end{conj}

We can refine Theorem  \ref{thm:inv-factorizations} to count the number of factorizations $\pi = \tau \circ \sigma$ where $\sigma$ and $\tau$ are involutions with $s$ and $t$ fixed points, respectively.  This requires determining all the possible unlabeled graphs on $[n]$ with properly $2$-colored edges which can be labeled to give two involutions which compose to a permutation with the cycle type of $\pi$, and then counting the labellings which actually give $\pi$.  This is impractical for large $s$ and $t$.  The fixed-point-free case, though, is straightforward.

\begin{prop} Let $c_1, c_2, \ldots, c_n$ be nonnegative even integers with $\sum_{k=1}^n kc_k = n$.  Then the number of factorizations of a permutation $\pi$ of type $1^{c_1} \ldots n^{c_n}$ into two fixed-point-free involutions is
\[ \prod_{k=1}^n (c_k-1)!! k^{c_k/2} \]
where we adopt the convention $(-1)!! = 1$. \end{prop}
\begin{proof} The graph $\sigma \cup \tau$ corresponding to such a factorization consists of $c_k$ cycles of length $2k$, for each $k$.  The permutation $k$-cycles can be paired up into graphical cycles in $(c_k-1)!!$ ways.  Each pair of permutation $k$-cycles thus obtained can be used to label a graphical cycle in any of $k$ ways, following Lemma \ref{lem:graphical-cycle}.  Thus the number of ways to arrange the elements of $k$-cycles of $\pi$ in the graphical representation is $(c_k-1)!! k^{c_k/2}$.  The total number of factorizations is just the product over cycle lengths. \end{proof}

We note that if any of the $c_k$ are odd, then $\pi$ has {\it no} factorizations into fixed-point-free involutions.  Furthermore, the proportion of permutations of $[n]$ having all $c_k$ even (that is, an even number of cycles of each length) is $\Theta(n^{-2})$.  The details of this enumeration and the cycle structure of such permutations will be considered in \cite{lugo-further}.


\begin{thebibliography}{99}
\bibitem{arratia} Richard Arratia. On the Stanley-Wilf Conjecture for the Number of Permutations Avoiding a Given Pattern.  {\it Electronic Journal of Combinatorics} 6 (1999) N1.
\bibitem{benaych-georges} Florent Benaych-Georges. Cycles of random permutations with restricted cycle lengths.   Preprint, arXiv:0712.1903.
\bibitem{billingsley} Patrick Billingsley.  {\it Probability and measure}, 3rd edition.  Wiley, 1995.
\bibitem{bona} Miklos Bona. The limit of a Stanley-Wilf sequence is not always rational and layered patterns beat monotone patterns. {\it Journal of Combinatorial Theory, Series A}, 110 (2005), 223-235.
\bibitem{bousquet-melou} Mireille Bousquet-Melou. Four Classes of Pattern-Avoiding Permutations Under One Roof: Generating Trees with Two Labels.  {\it Electronic Journal of Combinatorics} 9(2) (2003) R19.
\bibitem{dfls} P. Duchon, P. Flajolet, G. Louchard and G. Schaeffer. Boltzmann samplers for the random generation of combinatorial structures. {\it Combin. Probab. Comput.} 13 (2004), 577-625.
\bibitem{egge} Erik Egge.  Restricted Symmetric Permutations.  {\it Annals of Combinatorics} 11 (2007) 405-434.
\bibitem{egge-mansour} Erik Egge and Toufik Mansour.  231-avoiding involutions and Fibonacci numbers.  {\it Australasian Journal of Combinatorics} 30 (2004) 75-84.
\bibitem{flajolet-sedgewick} Philippe Flajolet and Robert Sedgewick.  {\it Analytic combinatorics.}  Cambridge, 2009.
\bibitem{gessel}  I. M. Gessel.  Symmetric functions and P-recursiveness. {\it Journal of Combinatorial Theory, Series A} 53 (1990) 257-285. 
\bibitem{greenhill-janson-rucinski} Catherine Greenhill, Svante Janson, Andrzej Rucinski.  On the number of perfect matchings in random lifts.  Preprint, arXiv:0907.0958.
\bibitem{goupil-schaeffer} Alain Goupil and Gilles Schaeffer.  Factoring n-cycles and counting maps of given genus. {\it European Journal of Combinatorics} 19 (1998) 819-834.
\bibitem{hayman} Walter Hayman.  A generalisation of Stirling's formula.  {\it Journal f\"ur die reine und angewandte Mathematik} 196 (1956) 67-95.
\bibitem{lugo} Michael Lugo.  Profiles of permutations.  {\it Electronic Journal of Combinatorics} 16(1) (2009) R99.
\bibitem{lugo-further} Michael Lugo.  Further examples of weighted permutations.  In preparation. 
\bibitem{marcus-tardos} A. Marcus, G. Tardos.  Excluded permutations, matrices and the Stanley-Wilf conjecture. {\it Journal of Combinatorial Theory Series A} 107 (2004) 153-160.
\bibitem{moser-wyman} L. Moser, M. Wyman.  On solutions of $x^d$ in symmetric groups.  {\it Canadian Journal of Mathematics} 7 (1955) 159-168.
\bibitem{pemantle-notes} Robin Pemantle, lecture notes for Math 581: Analytic combinatorics in more than one variable.  University of Pennsylvania, spring 2009.  Available online at \url{http://www.math.upenn.edu/\~pemantle/581-html/lecture-notes.html}.
\bibitem{regev} A. Regev.  Asymptotic values for degrees associated with strips of Young diagrams.  {\it Adv. Math.} 41:115-136, 1981. 
\bibitem{roberts-vivaldi} John A. G. Roberts, Franco Vivaldi.  A combinatorial model for reversible rational maps over finite fields.  Preprint, arXiv:0905.4135.
\bibitem{simion-schmidt} R. Simion and F. Schmidt. Restricted permutations.  {\it European Journal of Combinatorics} 6 (1985) 383-406.
\bibitem{stanley-99} Richard P. Stanley. {\it Enumerative Combinatorics}, volume 2.  Cambridge University Press, 1999.
\bibitem{timashev} A. N. Timashev. Random permutations with cycle lengths in a given finite set. {\it Discrete Mathematics and Applications} 18 (2008) 25-39.
\bibitem{west} Julian West.  Permutations with forbidden subsequences and stack-sortable permutations.  Ph. D. thesis, MIT, 1990.
\bibitem{wilf-bulletin} Herbert Wilf.  The asymptotics of $e^{P(z)}$ and the number of elements of each order in $S_n$.  {\it Bulletin of the American Mathematical Society} 15 (1986) 228-232.
\bibitem{wilf} Herbert Wilf. {\it generatingfunctionology}, 2nd edition.  Academic Press, 1994.  Available online at \url{http://www.math.upenn.edu/~wilf/DownldGF.html}.
\bibitem{wright-32} E. Maitland Wright.  The coefficients of a certain power series.  {\it J. London Math. Soc.} 7 (1932) 256-262.
\bibitem{wulcan} Elizabeth Wulcan.  Pattern avoidance in involutions.  Master's thesis, Chalmers University of Technology, 2002. Available online at \url{http://www.math.lsa.umich.edu/~wulcan/pattern.pdf}
\end{thebibliography}
\end{document}